\documentclass{article}
\usepackage{amsthm,amssymb,amsmath}
\usepackage{color}
\usepackage{epsfig,verbatim}

\newtheorem{thm}{Theorem}
\newtheorem{prop}[thm]{Proposition}

\newtheorem{conj}[thm]{Conjecture}

\theoremstyle{remark}

\newcommand{\dbar}{\overline{\partial}}

\title{The K\"ahler-Ricci flow and K-polystability}
\author{G\'abor Sz\'ekelyhidi}
\date{}

\begin{document}

\maketitle

\begin{abstract}
	We consider the K\"ahler-Ricci flow on a Fano manifold. We show
	that if the curvature remains uniformly bounded along the flow,
	the Mabuchi energy is bounded below, and the manifold is
	K-polystable, then the manifold admits a K\"ahler-Einstein
	metric. The main ingredient is a result that says that a
	sufficiently small perturbation of a cscK manifold admits a cscK
	metric if it is K-polystable.
\end{abstract}

\section{Introduction}
The question of existence of K\"ahler-Einstein metrics on Fano manifolds
is a fundamental problem in K\"ahler geometry. The case of K\"ahler
manifolds with zero or negative first Chern class was solved by
Yau~\cite{Yau78} and Aubin~\cite{Aub78} showing that every such manifold
admits a K\"ahler-Einstein metric. In the case of positive first Chern
class the presence of holomorphic vector fields led to the first obstructions
(Matsushima \cite{Mat57} and Futaki \cite{Fut83}). The problem can also be
generalised to manifolds $M$ with an ample line bundle $L$, where we can ask
whether there exists a constant scalar curvature (cscK) metric in the
class $c_1(L)$. The Futaki invariant generalises to this context (see
Calabi~\cite{Cal85}) and
since it enters later, we take this moment to define it. Let us choose a
K\"ahler metric $\omega$ on $M$ in the class $c_1(L)$, and write
$S(\omega)$ for the scalar curvature. Let the function $h$ satisfy
\[ S(\omega) - \hat{S} = \Delta_\omega h, \]
where $\hat{S}$ is the average of $S(\omega)$ with respect to the volume
form $\omega^n$. 
Given a holomorphic vector field $v$ on $M$ the Futaki invariant is
defined to be
\[ F(v) = \int_M v(h) \omega^n.\]
The main point is that $F(v)$ is independent of the metric we chose in
$c_1(L)$, and if there exists a cscK metric in $c_1(L)$ then clearly
$F(v)=0$ for all holomorphic vector fields $v$.

Later, inspired by a conjecture due to Yau \cite{Yau93}, a more subtle
obstruction called K-polystability has been found which is
related to stability in the sense of geometric invariant
theory. There are different notions of K-polystability, the original one
in the context of K\"ahler-Einstein manifolds due to Tian
\cite{Tian97}, and a more algebraic one which applies to the more general
problem of cscK metrics due to Donaldson
\cite{Don02}. The main ingredient in the definition is the notion of a
test-configuration for a polarised variety $(M,L)$. 
This is a $\mathbf{C}^*$-equivariant flat family
$(\mathcal{M},\mathcal{L})$ over $\mathbf{C}$, such that the generic
fibre is isomorphic to $(M,rL)$ for some $r>0$. The central fibre
$(M_0, L_0)$, which in Donaldson's definition might be an arbitrary
scheme, inherits a $\mathbf{C}^*$-action. This allows one to define the
Donaldson-Futaki invariant of the test-configuration, which coincides
with the Futaki invariant of the vector field generating the
$\mathbf{C}^*$-action when $M_0$ is smooth. We will only be
interested in the case when the central fibre is smooth, so we will not
give the detailed definition of the Donaldson-Futaki invariant.
The polarised manifold $(M,L)$ is called K-polystable if the
Donaldson-Futaki invariant is non-negative for all 
test-configurations and is zero
only for product configurations.
A product configuration is a test-configuration where the central
fibre is isomorphic to the general fibre. 

The main conjecture is the following. 
\begin{conj}[Yau-Tian-Donaldson] \label{conj:YTD}
	A polarised manifold $(M,L)$ admits a
	cscK metric in $c_1(L)$ if and only if it is K-polystable.
\end{conj}

In one direction Donaldson has shown that if $(M,L)$ admits a cscK
metric, then it is K-semistable. This is a weakening of K-polystability
where the Donaldson-Futaki invariant is allowed to be zero as well.
This was improved by
Stoppa~\cite{Sto08} who shows that existence of a cscK metric implies
K-polystability (using Donaldson's definition) 
when the manifold has no holomorphic vector fields. 
 In
the converse direction little progress has been made, except in the case
of toric surfaces, where Donaldson~\cite{Don08_1}
has proved Conjecture~\ref{conj:YTD}. In
the K\"ahler-Einstein case Tian proved that the existence of a
K\"ahler-Einstein metric implies his version of K-stability, where the
central fibre of test-configurations is only allowed to have normal
singularities. In the converse direction he showed that properness of a
certain energy functional (the Mabuchi functional) implies the existence
of a K\"ahler-Einstein metric, but it is still to be seen whether this
can be related to an algebraic condition. In the case of toric Fano
manifolds Conjecture~\ref{conj:YTD} is also known thanks to the work of
Wang-Zhu~\cite{WZ04}. Our first result is the
following.

\begin{thm}\label{thm:deform} Suppose $(M,L)$ is cscK, and let
	$(M',L')$ be a
	sufficiently small deformation of the complex structures of $M$ and
	$L$. If $(M',L')$ is K-polystable then it admits a cscK metric.
\end{thm}

It will be clear from the proof that we only need a weak version of
K-polystability where we only consider test-configurations with smooth
central fibres which are themselves small deformations of $(M,L)$. We
have been informed that T. Br\"onnle has obtained some similar results on
perturbing cscK metrics for his PhD thesis. When we are dealing with
Fano varieties polarised by the anticanonical bundle, then we do not
need to keep track of the polarization since it is canonically defined.
In particular the proof of the 
theorem shows that if $M$ is a Fano K\"ahler-Einstein
manifold, and $M'$ is a sufficiently small deformation of $M$, then
either $M'$ admits a K\"ahler-Einstein metric, or there is a
test-configuration for $M'$ with smooth central fibre $M''$. Moreover
$M''$ admits a K\"ahler-Einstein metric and it is itself a small
deformation of $M$ (or it can also be equal to $M$). 

We apply this result to the K\"ahler-Ricci flow on a Fano manifold. 
It was shown by Cao~\cite{Cao85}
that in the case of negative or zero first Chern class the flow
converges to the K\"ahler-Einstein metric which is guaranteed to exist
by Yau and Aubin's theorems.
In the Fano case Cao showed that the flow exists for all time, and the
main problem is to find conditions under which it converges. In view of
the Yau-Tian-Donaldson conjecture one would like to show that under the
assumption of K-polystability the flow converges to a K\"ahler-Einstein
metric. Without additional assumptions this seems out of reach at
present. One interesting question is what we can say if we assume that
the Riemannian curvature is uniformly bounded along the flow. The main
result about the K\"ahler-Ricci flow that we use is the following, which 
is based on Perelman's work~\cite{Per02}. 

\begin{thm}[see~\cite{PSSW07_2}, \cite{PS06_1}]\label{thm:converge} 
	Suppose that the Riemann 
	curvature tensor is uniformly bounded
	along the K\"ahler-Ricci flow on a Fano manifold $M$. Write $J$
	for the complex structure of $M$. We can
	then find a sequence of diffeomorphisms $\phi_k:M\to M$ such
	that $\phi_k^*(J)$ converges in $C^\infty$ to a smooth complex
	structure $J_0$, such that $(M,J_0)$ admits a K\"ahler-Ricci
	soliton. If in addition the Mabuchi functional of $M$ is bounded
	below, then $(M,J_0)$ admits a K\"ahler-Einstein metric.
\end{thm}

The question is then to find conditions which ensure that the complex
structure $J_0$ is isomorphic to $J$.
Phong and Sturm \cite{PS06_1} have introduced such a condition
called Condition B. 
Let us briefly recall that a complex manifold $(M,J)$ satisfies
Condition B if we cannot find a sequence of diffeomorphisms $\phi_k$
with $\phi_k^*(J)$ converging in $C^\infty$ to a complex structure $J_0$
which has a strictly higher dimensional space of
holomorphic vector fields than $J$. 
The result in \cite{PS06_1} 
was improved in \cite{PSSW07_2} and it is the following.

\begin{thm}[Phong-Song-Sturm-Weinkove] Suppose that the Riemann
	curvature tensor is uniformly bounded along the K\"ahler-Ricci
	flow on a Fano manifold $M$. If the Futaki invariant of $M$
	vanishes, and $M$ satisfies Condition B, then the flow converges
	exponentially fast in $C^\infty$ to a K\"ahler-Einstein metric.
\end{thm}

In this
direction we will prove the following related result. 

\begin{thm}\label{thm:condB}
	If the Riemann curvature tensor is uniformly bounded along the
	K\"ahler-Ricci flow on a Fano manifold $M$, and $M$ satisfies
	Condition B, then $M$ admits a K\"ahler-Ricci soliton.
\end{thm}

The proof of this does not need Theorem~\ref{thm:deform}.
Note that if the Futaki invariant of $M$ with respect to the class
$c_1(M)$ vanishes, then a K\"ahler-Ricci
soliton on $M$ is necessarily a K\"ahler-Einstein metric. 
Our main result is that under the additional assumption that the Mabuchi
functional is bounded from below, we can replace Condition B with
K-polystability.

\begin{thm}\label{thm:KRflow}
	Suppose that the Riemann curvature tensor is uniformly
	bounded along the K\"ahler-Ricci flow on a Fano manifold $M$.
	Suppose in addition that the Mabuchi functional on $M$ is
	bounded from below and that $M$ is K-polystable. Then $M$
	admits a K\"ahler-Einstein metric.
\end{thm}

In the Fano case when 
saying that $M$ is K-polystable we mean that the pair $(M,-K_M)$
is K-polystable. 
Note that once the existence of a K\"ahler-Einstein metric is
established, convergence of the flow follows from the work of Tian and
Zhu~\cite{TZ07}. Combined with the result of
Phong-Song-Sturm-Weinkove~\cite{PSSW07} we obtain exponential
convergence. 
We hope that the assumption that the Mabuchi functional is bounded from
below can be removed in the future. Note that for toric varieties
Donaldson~\cite{Don02} has shown that K-polystability implies that the
Mabuchi functional is bounded from below. 

We will see that instead of $M$ being K-polystable it is enough to
require that there 
are no non-product test-configurations for $M$ with smooth central fibre
which have zero Futaki invariant. The lower bound for the Mabuchi
functional should be thought of as a semistability condition, which is
strengthened to stability by excluding test-configurations with zero
Futaki invariant. Bounding
the curvature along the flow allows us to 
consider test-configurations
with smooth central fibres. Also, the main difficulty in
replacing Condition B with K-polystability is to obtain instead of just a
sequence of diffeomorphisms $\phi_k$ with $\phi_k^*(J)\to J_0$, a
test-configuration with central fibre $J_0$. 

In the next section we give the proof of Theorem \ref{thm:condB}, and
assuming Theorem \ref{thm:deform} the proof of Theorem~\ref{thm:KRflow}.
The rest of the
paper will then be devoted to the proof of Theorem~\ref{thm:deform}.

\subsubsection*{Acknowledgements}
I would like to thank J. Stoppa and V. Tosatti for many 
helpful discussions and X. X. Chen for pointing out a mistake in an
earlier version of the paper.  
The ideas in this paper have much to owe to S. K. Donaldson's
view of the subject and I have benefitted greatly from conversations with
him. In addition I am grateful to D. H. Phong for his interest in my
work and his encouragement. This work was carried out while visiting
Harvard University and was supported by an EPSRC Postdoctoral
Fellowship at Imperial College London. 

\section{Proofs of Theorems \ref{thm:condB} and \ref{thm:KRflow}}

\begin{proof}[Proof of Theorem~\ref{thm:condB}]
	Under our assumptions Theorem~\ref{thm:converge} implies that we
	can find a sequence of diffeomorphisms $\phi_k$, such that that
	sequence of complex structures $J_k=\phi_k^*(J)$ converges to a
	complex structure $J_0$, and $(M,J_0)$ admits a K\"ahler-Ricci
	soliton. We want to show that Condition B implies that $J_0$ is
	isomorphic to $J$. 

	We have $c_1(M,J_k)
	\to c_1(M, J_0)$ and since these are integral cohomology
	classes we must have $c_1(M,J_k) = c_1(M, J_0)$ for
	sufficiently large $k$.
	This means that we can think of the canonical
	bundles $K_{J_k}$ of $(M,J_k)$ as a fixed complex line bundle over $M$
	with complex structures varying with $k$. We can then take a
	basis of sections of $-lK_{J_0}$ for some large $l$ which gives
	an embedding of $(M,J_0)$ into projective space, and perturb it
	to a basis of sections of $-lK_{J_k}$ for large $k$. This will
	give a sequence of embeddings $V_k\subset\mathbf{P}^N$ 
	of $(M,J)$, which converges to an embedding
	$V\subset\mathbf{P}^N$ of $(M,J_0)$. 

	The fact that all the $V_k$ are isomorphic to each other implies
	that they are all in the same orbit of $PGL(N+1,\mathbf{C})$
	acting on the Hilbert scheme. The fact that the $V_k$ converge
	to $V$ means that $V$ represents a point in the closure of this
	orbit. If $J_0$ is not isomorphic to $J$, then $V$ represents a
	point in the boundary of the orbit of $V_k$, in which case its
	stabiliser subgroup must have dimension strictly greater than
	that of $V_k$. This is because the boundary of an orbit is a
	union of strictly lower dimensional orbits. The stabiliser in
	this case is just the group of holomorphic automorphisms, so
	this contradicts the assumption that $J$ satisfies Condition B.
	Therefore $J_0$ is isomorphic to $J$, and so $(M,J)$ admits a
	K\"ahler-Ricci soliton. 
\end{proof}

Assuming Theorem \ref{thm:deform} we can give the proof of Theorem
\ref{thm:KRflow}.

\begin{proof}[Proof of Theorem \ref{thm:KRflow}]
	In this case Theorem~\ref{thm:converge} implies that
	there exists a sequence of diffeomorphisms $\phi_k$ such that
	$J_k=\phi_k^*(J)$ converges to a complex structure
	$J_0$ which this time admits a K\"ahler-Einstein metric. As in
	the previous proof we have $c_1(M,J_k)=c_1(M,J_0)$
	for sufficiently large $k$. This means that we can apply
	Theorem~\ref{thm:deform}. Since by assumption $(M,J_k)$ is
	K-polystable, the theorem implies that $(M,J_k)$ admits a
	cscK metric in the class $c_1(M,J_k)$ for sufficiently large $k$. 
	This is necessarily a K\"ahler-Einstein metric on $(M,J)$. 
\end{proof}

\section{Perturbing cscK metrics}
Suppose that the complex manifold $(M,J)$ admits a cscK metric $\omega$.
In this section we study the problem of whether small deformations $J_t$
of the complex structure admit a cscK metric. 
We restrict attention to deformations of the complex structure which are
compatible with $\omega$. When $(M,J)$ admits no holomorphic vector fields,
then the implicit function theorem
shows that every small such deformation admits a cscK
metric. The case when $(M,J)$ has holomorphic vector fields is more
subtle, and was also studied by LeBrun-Simanca~\cite{LS94} in relation
with the Futaki invariant. 

Since cscK metrics can be interpreted
as zeros of a moment map, the deformation
problem can be cast into a general framework.  We first
recall this moment map picture from Donaldson \cite{Don97}.

\subsubsection*{The moment map picture}
Let us write $\mathcal{J}$ for the space of almost complex structures on
$M$, compatible with $\omega$. The tangent space $T_J\mathcal{J}$
at a point $J$ can be
identified with the space of 1-forms
$\alpha\in\Omega^{0,1}(T^{1,0})$ which satisfy
\[ \omega(\alpha(X),Y) + \omega(X,\alpha(Y))=0.\]
This space has a natural complex structure and also an $L^2$ inner
product, which gives $\mathcal{J}$ the structure of an infinite
dimensional K\"ahler manifold. Let us write $\mathcal{G}$ for the group
of exact symplectomorphisms of $(M,\omega)$. This acts naturally on
$\mathcal{J}$ preserving the K\"ahler structure. The Lie algebra of
$\mathcal{G}$ can be identified with $C^\infty_0(M)$ via the Hamiltonian
construction, and we identify it with its dual using the $L^2$ product
induced by $\omega$.
It was shown by
Donaldson (also Fujiki \cite{Fuj92}) that the map
\begin{equation}\label{eq:mmap}
	J \mapsto S(J,\omega) - \hat{S} 
\end{equation}
is an equivariant moment map for this action. Here $S(J,\omega)$ is the
``Hermitian scalar curvature'' defined in \cite{Don97}, which when
$J$ is integrable coincides with the usual scalar curvature of the
K\"ahler metric defined by $(J,\omega)$ up to a
constant factor. Moreover $\hat{S}$ is the average of $S(J,\omega)$,
which is independent of $J$. We see therefore that if $J$ is integrable
and is a zero of this moment map, then $(J,\omega)$ defines a cscK
metric.

To explain what is meant by (\ref{eq:mmap})
being a moment map, define the following two operators at
$J\in\mathcal{J}$.
The infinitesimal action of $C^\infty_0(M)$ is given by
\[ P : C^\infty_0(M) \to T_J\mathcal{J}, \]
which we can also write as $P(H) = \dbar X_H$ where $X_H$ is the
Hamiltonian vector field corresponding to $H$. 
The other operator is the derivative of $S(J,\omega)$, 
\[ Q : T_J\mathcal{J} \to C^\infty_0(M).\]
The fact that (\ref{eq:mmap}) gives a moment map can be expressed as
\begin{equation}\label{eq:mmap1}
	\langle Q(\alpha), H \rangle_{L^2} = \Omega(\alpha, P(H)),
\end{equation}
where $\Omega$ is the symplectic form on $\mathcal{J}$. 

\subsubsection*{The complex orbits}
A key observation in \cite{Don97} is that while the complexification
$\mathcal{G}^c$ of the group $\mathcal{G}$ does not exist, one can
still make sense of its orbits as leaves of a foliation, and a leaf
containing an integrable complex structure can be interpreted as
the space of K\"ahler metrics in a K\"ahler class. 
We briefly explain how this works. 
We can complexify the action of $\mathcal{G}$ on $\mathcal{J}$
on the level of Lie algebras
by extending the operator $P$ to 
\[ P : C^\infty_0(M,\mathbf{C}) \to T_J\mathcal{J}\]
in the natural way. We can
then think of leaves of the resulting foliation on $\mathcal{J}$ as
the orbits of $\mathcal{G}^c$. With this in mind
we will say that $J_0$ and $J_1$ are in the same $\mathcal{G}^c$-orbit if we
can find $\phi_t\in C^\infty_0(M,\mathbf{C})$ and a path $J_t\in\mathcal{J}$
for $t\in[0,1]$ joining $J_0$ and $J_1$, which satisfies
\[ \frac{d}{dt} J_t = P_t(\phi_t).\]
We write $P_t$ to emphasise that the operator $P$ depends
on the complex structure. When $J_0$ is
integrable, then in fact there exists a diffeomorphism $f:M\to M$ and
some $\psi\in C^\infty_0(M)$ such
that $f^*(J_1) = J_0$ and $f^*(\omega)=\omega+i\partial\dbar\psi$. This
means that up to the action of diffeomorphisms, integrable complex
structures in the same $\mathcal{G}^c$ orbit can be thought of as
K\"ahler metrics in the same K\"ahler class on a fixed complex manifold.

We will later need to perturb complex structures in a $\mathcal{G}^c$
orbit of an integrable $J$, 
so we give the relevant definition here. Let $J$ be an almost complex
structure compatible with $\omega$ and let
$U\subset L_k^2$ be a small ball around the origin for some sufficiently
large $k$. For $\phi\in U$
we can define a complex structure $F_\phi(J)$ in the following way.
For $t\in[0,1]$ write
\[ \omega_t = \omega - tdJd\phi.\]
Then 
\[ \frac{d}{dt} \omega_t = d\alpha\]
for a fixed one-form $\alpha$ with coefficients in $L_{k-1}^2$. 
Write $X_t$ for the vector field dual to
$-\alpha$ under the symplectic form $\omega_t$ (so $X_t$ has
coefficients in $L_{k-2}^2$). Then
\[ \frac{d}{dt} \omega_t = -d(\iota_{X_t}\omega_t) = -L_{X_t}\omega_t.\]
We can now define diffeomorphisms $f_t$ with coefficients in $L_{k-2}^2$ 
by integrating the family of
vector fields $X_t$ for $t\in[0,1]$. 
\[ \frac{d}{dt} f_t = X_t. \]
Then $f_1^*(\omega_1)=\omega$, and we let $F_\phi(J) =f_1^*J$ which has
coefficients in
 $L_{k-2}^2$. 
Note that when $J$ is integrable then the two K\"ahler manifolds
$(J,\omega - dJd\phi)$ and $(F_\phi(J), \omega)$ are
isometric, so up to a diffeomorphism, we are just perturbing the metric
in its K\"ahler class.

The map $\phi \mapsto F_\phi(J)$ obtained this way is $K$-equivariant,
where $K$ is the stabiliser of $J$ in $\mathcal{G}$. Here a
diffeomorphism in $K$ acts on
$U\in L_k^2$ linearly by pulling back the functions. 

We will later need to know the derivative of the map $\phi\mapsto
F_\phi(J)$ at the
origin. This is a map from $U\in L_k^2$ to complex structures in
$L_{k-2}^2$. By
the construction the derivative at the origin is given by
\begin{equation} \label{eq:DF}
	DF_0(\phi) = J(L_{X_\phi} J) = JP(\phi),
\end{equation}
where $X_\phi$ is the Hamiltonian vector field corresponding to $\phi$.

\subsubsection*{Construction of a local slice}
Suppose now that $(J_0,\omega)$ is a cscK metric on $M$ (in particular
$J_0$ is integrable). 
Following Kuranishi \cite{Kur65} we can construct a local
slice for the action of $\mathcal{G}^c$ on $\mathcal{J}$ near $J_0$, which
intersects the $\mathcal{G}^c$ orbit of every integrable
complex structure near $J_0$. We will also allow some non-integrable complex
structures in the slice and this will allow the slice to be smooth.

Recall that the infinitesimal
action of $\mathcal{G}^c$ is given by the complexification of $P$,
\[ P : C^\infty_0(M,\mathbf{C}) \to T_{J_0}\mathcal{J}.\]
We also have an operator $\dbar: T_{J_0}\mathcal{J}
\to \Omega^{0,2}(T^{1,0})$, and the two fit into an elliptic complex
\begin{equation}\label{eq:ell}
	C^\infty_0(M,\mathbf{C}) \stackrel{P}{\longrightarrow}
	T_{J_0}\mathcal{J} \stackrel{\dbar}{\longrightarrow}
\Omega^{0,2}(T^{1,0}). 
\end{equation}
Let us write 
\[ \tilde{H}^1 = \{ \alpha\in T_{J_0}\mathcal{J}\,|\, P^*\alpha=\dbar\alpha=0\}. \]
This is a finite dimensional vector space since it is the kernel of the
elliptic operator $PP^*+(\dbar^*\dbar)^2$ on $T_{J_0}\mathcal{J}$ 
(for more details
see~\cite{FS90}). We will write $K$ for the
stabiliser of $J_0$ in $\mathcal{G}$, ie. the group of Hamiltonian
isometries of $(J_0,\omega)$, and $\mathfrak{k}$ for its Lie algebra.
Note that $\mathfrak{k}$ can be identified with the kernel of $P$ in
$C^{\infty}_0(M,\mathbf{R})$. 
The group $K$ acts naturally on $\tilde{H}^1$, and we write $K^c$ for
the complexification of $K$. 

\begin{prop}\label{prop:slice}
	There exists a ball $B\subset \tilde{H}^1$ around the origin
	and a
	$K$-equivariant map 
	\[ \Phi : B \to \mathcal{J}, \]
	such that the $\mathcal{G}^c$ orbit of every integrable $J$ near
	$J_0$ intersects the image of $\Phi$. Also
	if $x,x'$ are in the same
	$K^c$-orbit and $\Phi(x)$ is an integrable complex structure,
	then $\Phi(x), \Phi(x')$ are in the same
	$\mathcal{G}^c$-orbit. Moreover for all $x\in B$ we have
	 $S(\Phi(x),\omega)\in\mathfrak{k}$.
\end{prop}
We learned the idea of requiring this 
last condition in order to reduce the problem
to a finite dimensional one from \cite{Don08}.
\begin{proof}
	Following Kuranishi~\cite{Kur65}
	we can construct a $K$-equivariant holomorphic map 
\[ \Phi_1 : B_1 \to \mathcal{J}\]
from some ball $B_1$ in $\tilde{H}^1$, such that the
$\mathcal{G}^c$-orbit of every integrable complex structure near $J_0$
intersects the image of $\Phi_1$. The difference is that in Kuranishi's
situation instead of (\ref{eq:ell}) the relevant elliptic complex is
\[ \Gamma(T^{1,0}) \stackrel{\dbar}{\longrightarrow}
\Omega^{0,1}(T^{1,0}) \stackrel{\dbar}{\longrightarrow} 
\Omega^{0,2}(T^{1,0}),
\]
since he is constructing a slice for the action of the diffeomorphism
group instead of $\mathcal{G}^c$. In addition we do not insist that
all the complex structures in the image of
$\Phi$ should be integrable; the integrable ones will correspond to an
analytic subset of $B$ in the same way as in \cite{Kur65}.

It only remains to show that we can perturb $\Phi_1$ in
such a way as to satisfy the last statement in the proposition. We will
perturb inside the $\mathcal{G}^c$ orbits using the map $F_\phi$ defined
in the previous subsection. For
$J\in\mathcal{J}$ sufficiently close to $J_0$ and $\phi\in L_l^2$
sufficiently small, we have a complex structure $F_\phi(J)$ in
$L_{l-2}^2$ if $l$ is large enough.

Let
us write $\mathfrak{k}^\perp_l$ for the orthogonal complement to
$\mathfrak{k}$ in the Sobolev space $L^2_l$, let $U_l\subset
\mathfrak{k}_l^\perp$ be a
small ball around the origin and consider the map
\[ \begin{aligned}
	G : B_1\times U_l &\to \mathfrak{k}^\perp_{l-4}
	\\ 
	(x,\phi) &\mapsto \Pi_{\mathfrak{k}^\perp_{l-4}}
	S(F_\phi(\Phi_1(x)), \omega),
\end{aligned} \]
where $\Pi$ is the $L^2$ orthogonal projection. 
It follows from (\ref{eq:mmap1}) and (\ref{eq:DF}) that 
the derivative of $G$ at the origin is given by 
\[ DG_{(0,0)}(\phi) = P^*P(\phi).\]
This is an isomorphism from $\mathfrak{k}_l^\perp\to\mathfrak{k}_{l-4}^\perp$,
so by the implicit function theorem we can perturb $\Phi_1$ to 
\[ \Phi : B \to\mathcal{J},\]
where $B$ is a smaller ball than $B_1$, and so that for all $x\in B$ we
have $S(\Phi(x),\omega)\in\mathfrak{k}$.  Moreover $F$ and $\Phi_1$ are
$K$-equivariant, therefore so are $G$ and $\Phi$.
\end{proof}

Let us write $\Omega$ for the symplectic form on $B$ pulled back from
$\mathcal{J}$ via $\Phi$. This form is preserved by the $K$-action on
$B$, and a moment map for the action is given by
\[ \mu(x) = S(\Phi(x),\omega)\in\mathfrak{k}.\]
Moreover points $x,x'$ in the same $K^c$ orbit correspond to 
complex structures in the same $\mathcal{G}^c$-orbit if they represent
integrable complex structures. Also note that the
$K$ and $K^c$ actions on $B$ are just the linear ones induced by those
on $\tilde{H}^1$. We have therefore reduced our problem to finding
$K^c$-orbits which contain zeros of the moment map $\mu$. 

\subsubsection*{The finite dimensional problem}

We want to prove the following.
\begin{prop}\label{prop:finitedim}
	After possibly shrinking $B$, suppose that $v\in B$ is
	polystable for the $K^c$-action on $\tilde{H}^1$. Then there is
	a $v_0\in B$ in the $K^c$-orbit of $v$ such that $\mu(v_0)=0$. 
\end{prop}
\begin{proof}
	Let us identify $\tilde{H}^1$ with the tangent space to $B$ at
	the origin, and write $\Omega_0$ for the linear symplectic form
	induced on $\tilde{H}^1$ by $\Omega$. Also write 
	\[ \nu : \tilde{H}^1 \to \mathfrak{k} \]
	for the corresponding moment map, where we have identified
	$\mathfrak{k}$ with its dual using the inner product induced by
	the $L^2$ product on functions as before. 

	If $v\in B$ is polystable for the $K^c$-action then by the
	Kempf-Ness theorem there is a zero of the moment map $\nu$ in
	the $K^c$-orbit of $v$. In fact this is obtained by minimising
	the norm over the $K^c$ orbit, so the zero of the moment map
	will still be in $B$. We can therefore assume without loss of
	generality that $\nu(v)=0$. 

	We have the Taylor expansion
	\[\mu(tv) = \mu(0) + t\,d\mu_0(v) +
	\frac{t^2}{2}\left.\frac{d^2}{dt^2}\right|_{t=0}\mu(tv) 
	+ O(t^3).\]
	By assumption $\mu(0)=0$ and since $0$ is a fixed point of the
	$K$-action we have $d\mu_0(v)=0$. Also it is easy to check
	that
	\[ \left.\frac{d^2}{dt^2}\right|_{t=0}\mu(tv) = \nu(v),\]
	so in sum we have $\mu(tv) = O(t^3)$. 

	For $x\in B$ let us write $K_x$ for the
	stabiliser of $x$ and $\mathfrak{k}_x$ for its Lie algebra. Note
	that $K_{tx}=K_x$ for all nonzero $t$. Then
	for all $\xi\in\mathfrak{k}_x$ we have
	\[ \frac{d}{dt}\langle\mu(tx),\xi\rangle = \Omega_{tx}(x,
	\sigma_{tx}(\xi)) = 0, \]
	where $\sigma_x : \mathfrak{k}\to T_xB$ is the infinitesimal
	action and we have identified the tangent space to $B$ at every
	point with $\tilde{H}^1$. This implies that for all $x\in B$ we have 
	\begin{equation} \label{eq:fut}
		\mu(x)\in \mathfrak{k}_x^\perp.
	\end{equation}

	We now try to perturb $tv$ into a zero of $\mu$ in its $K^c$
	orbit using Proposition~\ref{prop:pert} below for sufficiently small
	$t$. Write $Q_x$ for the operator
	\[\sigma_x^*\sigma_x : \mathfrak{k}_x^\perp \to
	\mathfrak{k}_x^\perp.\] 
	We need to estimate the norm $\Vert Q_x^{-1}\Vert$ for
	$x=e^{i\xi}\cdot(tv)$ where $\Vert\xi\Vert < \delta$ for some
	fixed small $\delta>0$. Clearly
	there is a constant $C$ such that 
	\[ \Vert\sigma_x(\eta)\Vert_{\Omega_0}^2 \geqslant C\Vert\eta\Vert^2,
	\quad\text{for all }\eta\in\mathfrak{k}_x^\perp\text{ and }
	x=e^{i\xi}v,\text{ where }\Vert\xi\Vert < \delta.\]
	Note that we have used the metric induced by $\Omega_0$.
	But $\sigma_{tx}(\eta)=t\sigma_x(\eta)$, and if the ball $B$ is
	sufficiently small then the metric on $B$ induced by $\Omega$
	is bounded below by half of the metric induced by
	$\Omega_0$, so we have
	\[ \Vert \sigma_{tx}(\eta)\Vert_{\Omega} \geqslant \frac{1}{2}
	t\Vert\sigma_x(\eta)\Vert_{\Omega_0}.\]
	It follows that for all $x=e^{i\xi}v$ with $\Vert\xi\Vert<\delta$ we
	have
	\[ \langle\sigma_{tx}^*\sigma_{tx}(\eta),\eta\rangle =
	\Vert\sigma_{tx}(\eta)\Vert^2_{\Omega} \geqslant
	\frac{1}{4} t^2 C\Vert\eta\Vert^2,\]
	and so $\Vert Q_{tx}^{-1}\Vert < C_1t^{-2}$.

	Now Proposition~\ref{prop:pert} implies that if
	$C_1t^{-2}\Vert\mu(tv)\Vert
	< \delta$, then there is a zero of $\mu$ in the $K^c$-orbit of
	$tv$. Since $\mu(tv) = O(t^3)$, this will be true for
	sufficiently small $t$, and so the proof of the proposition is
	complete.
\end{proof}

We have used the following extension of a result in \cite{Don01}. As in
the previous proof, let us write $Q_x$ for the operator
$\sigma_x^*\sigma_x$. This is an isomorphism 
\[ Q_x : \mathfrak{k}_x^\perp \to \mathfrak{k}_x^\perp,\]
and we write $\Lambda_x$ for the norm of $Q_x^{-1}$. 
\begin{prop}\label{prop:pert}
	Suppose $x_0\in B$ satisfies
	$\mu(x_0)\in\mathfrak{k}_{x_0}^\perp$. Given real numbers
	$\lambda,\delta$ such that $\Lambda_x\leqslant\lambda$ for all
	$x=e^{i\xi}x_0$ with $\Vert\xi\Vert < \delta$, suppose that
	$\lambda\Vert\mu(x_0)\Vert < \delta$. Then there is a point
	$y=e^{i\eta}x_0$ with $\mu(y)=0$, where $\Vert\eta\Vert\leqslant
	\lambda\Vert\mu(x_0)\Vert$.
\end{prop}

In Donaldson's statement it is assumed that $\mathfrak{k}_{x_0}$ is
trivial, but the proof of this slightly more general result is identical. 
We can now give the proof of Theorem~\ref{thm:deform}.

\begin{proof}[Proof of Theorem \ref{thm:deform}]
	Let $(J,\omega)$ be a cscK metric on $M$ in the class $c_1(L)$. 
	If $(M',L')$ is a small deformation of $(M,L)$ then necessarily
	$c_1(L)=c_1(L')$ as cohomology classes on the underlying smooth
	manifold. In particular $c_1(L)$ is a $(1,1)$-class
	with respect to the complex structure of $M'$ and by modifying
	$J'$ and $L'$ by a small diffeomorphism we can assume
	that $J'$ is compatible with $\omega$. 
	 Any sufficiently small
	deformation $J'$ of $J$ which is integrable and compatible with
	$\omega$
	is represented by some $v\in B$. If $v$ is polystable for the
	$K^c$-action, then Propositions~\ref{prop:slice}
	and~\ref{prop:finitedim} imply
	that $J'$ admits a cscK metric in the
	K\"ahler class $[\omega]$. What remains to be shown is that 
	if $v$ is not polystable then $(M',L')$ is not K-polystable. By the
	Hilbert-Mumford criterion there exists a one-parameter subgroup
	$\rho: \mathbf{C}^*\to K^c$ such that 
	\[ v_0 = \lim_{\lambda\to 0}\rho(\lambda)\cdot v \]
	is polystable (in fact it is a zero of the moment map for the
	linear symplectic form). Moreover $v_0\in B$, and it represents
	an integrable complex structure $J_0$ since integrability is a
	closed condition. Also, $J_0$ admits a cscK
	metric in $c_1(L)$ and so $(M,J_0)$ has vanishing Futaki
	invariant (with respect to the polarisation $c_1(L)$). 
	
	We can now construct a
	test-configuration for $(M,J')$ whose central fibre is
	$(M,J_0)$ as follows. First using $\Phi_1$ from the proof of
	Proposition~\ref{prop:slice} together with $\rho$, we obtain an
	$S^1$-equivariant holomorphic map 
	\[ F : \Delta \to\mathcal{J} \]
	from a small disk $\Delta$, such that $F(t)$ is isomorphic to
	$J'$ for non-zero $t$, and $F(0)$ is isomorphic to $J_0$. We let
	the total space of our test-configuration be
	$\mathcal{M}=M\times\Delta$ as
	a smooth manifold, endowed with the almost complex structure
	which on the fibre $M\times\{t\}$ is given by $F(t)$. Since $F$
	is holomorphic, this gives an integrable complex structure on
	$\mathcal{M}$. The $S^1$-action on $\mathcal{M}$ is the
	product action
	\[ \lambda\cdot(x,t) = (\rho(\lambda)\cdot x,\lambda t),\]
	where we have identified elements of $K$ with diffeomorphisms
	of $M$ fixing $J_0$ (we can assume that $\rho(S^1)\subset K$).
	Let us write $\overline{L}$ for the complex vector bundle
	underlying $L$ and $L'$. Let us also fix a connection $\nabla$
	on $\overline{L}$ so that the $(0,1)$-part of $\nabla$ with
	respect to the complex structure $J'$ gives $\overline{L}$ the
	holomorphic structure $L'$. Note that $\rho$ induces an
	$S^1$-action on $\overline{L}$ and we can assume that $\nabla$
	is $S^1$-invariant. We now let $\mathcal{L}$ be the pullback
	$\pi^*(\overline{L})$ under the projection $\pi:M\times\Delta\to
	M$, and endow it with the pullback connection and induced
	holomorphic structure (using the complex structure on
	$\mathcal{M}$ that we have defined). Then $\mathcal{L}$ has a
	natural $S^1$-action lifted from the action on $\mathcal{M}$,
	which preserves the connection and hence is
	holomorphic. The restriction of $\mathcal{L}$ to $M\times\{t\}$
	for non-zero $t$ is just $L'$ and while the restriction to
	$M\times\{0\}$ might not be $L$ as a holomorphic bundle, it at
	least has the same first Chern class and is therefore ample. We
	have thus obtained a flat, polarised, 
	$S^1$-equivariant family over $\Delta$ with general fibre
	$(M,J')$ polarised by $L'$, and central fibre $(M,J_0)$. The
	$S^1$-action extends to a $\mathbf{C}^*$-action, using which we
	can extend the family to a family over $\mathbf{C}$ so we have a
	test-configuration. Since 
	the Futaki invariant of any vector field on $(M,J_0)$ vanishes,
	this test-configuration has zero Futaki invariant. Also
	$(M,J_0)$ has more holomorphic vector fields (since the
	dimension of the stabiliser 
	$\mathfrak{k}_{v_0}$ is greater than that of
	$\mathfrak{k}_v$), so it is not
	isomorphic to $(M,J')$ and the test-configuration is not a
	product configuration. This shows that $(M,J')$ is not
	K-polystable. 
\end{proof}

\end{document}